\documentclass[11pt,reqno]{article}
\usepackage{amsmath,amsthm}

\input{tcilatex}

\begin{document}

\begin{center}
{\large K\"{a}hler Calculus: Quarks as Idempotents in}\linebreak {\large \
Solutions with Symmetry of Exterior Systems\bigskip }

Jos\'{e} G. Vargas

PST Associates

{138 Promontory Rd., Columbia, SC29209}\linebreak {josegvargas@earthlink.net%
\bigskip }

{To the memory of Erich K\"{a}hler \bigskip }
\end{center}

\textbf{Abstract.} Symmetry in differential equations constitutes a
significant starting point for theoretical arguments. With his Clifford
algebra of differential forms, K\"{a}hler addresses the overlooked
manifestation of symmetry in the solutions themselves of such equations. In
``K\"{a}hler's algebra'', solutions with a given symmetry are members of
left ideals generated by corresponding idempotents. These combine with phase
factors to take care of all the dependence of the solution on $x^{i}$ and $%
dx^{i}$ for each $\partial /\partial x^{i}=0$ symmetry. The maximum number
of idempotents (each for a one-parameter group) that can, therefore, go into
a solution is the dimensionality $n$ of the space. This number is further
limited by non-commutativity of idempotents.

We consider the tensor product of K\"{a}hler algebra with tangent Clifford
algebra, i.e. of valuedness. It has a commutative subalgebra of ``mirror
elements'', i.e. where the valuedness is dual to the differential form (like
in $dx\mathbf{i}$, but not in $dx\mathbf{j}$). It removes the aforementioned
limitation in the number of idempotent factors that can simultaneously go
into a solution. More importantly, it geometrizes the unit imaginary.

We carry K\"{a}hler's treatment over to the aforementioned commutative
algebra. Non-commutativity thus ceases to be a limiting factor in the
formation of products of idempotents representing each a one-parameter
symmetry. The interplay of factors that emerges in the expanded set of such
products is inimical to operator theory. For practical reasons, we stop at
products involving three one-parameter idempotents ---henceforth call
ternary--- even though the dimensionality of the all important spacetime
manifold allows for four of them.

\section{Introduction}

The present paper is an argument about manifestation of symmetry in
solutions of equations, not in the equations themselves. It results from
Erich K\"{a}hler's proposed treatment of symmetry of solutions of exterior
systems of differential equations. It goes beyond the traditional approach,
since phase factors are not enough to deal with this issue. Products of
idempotents are required, each of them associated with a corresponding
one-parameter group of symmetries. That is a deeper treatment of coexistence
of different symmetries in a solution than through the commutativity or not
of corresponding operators.

The term exterior system evokes the concept of exterior algebra. Although
the term idempotent presupposes a dot product, we have emphasized the dot
product in the title to compensate for the fact that the term exterior might
indicate its absence and thus cause confusion. As is well known, elements of
a Clifford algebra can always be written as sums of interior products. And
it is a historical accident that the name exterior system has been used
whenever a system of differential equations is written in terms of exterior
differential forms.

In section 2, we deal with Kaehler's treatment of translation and rotation
symmetry, specifically with their products. But, as relevant because of its
potential implications, is K\"{a}hler's study of angular momentum and of
total angular momentum on idempotents in his algebra of scalar-valued
differential forms. We have in mind the potential richness of theoretical
results when the same line of argument is carried out with Clifford valued
differential forms, which we introduce below. K\"{a}hler dealt at length not
only with scalar-valuedness but also with tensor valuedness, but it did not
have consequence for his treatment of solutions of exterior systems. We do
not deal with that topic here, since it would lead us directly into physics,
which we avoid in this paper.

In section 3, we consider symmetry related transformations from a
perspective of Clifford algebra of scalar-valued differential forms. When
the information they provide is taken to the idempotents, the unit imaginary
is geometrized. Actually there are many geometric units of square minus one,
which automatically increases the number of available idempotents without
resort to new symmetries. Tangent algebra valuedness thus makes a natural
entrance into the study of solutions with symmetry of exterior systems.

We keep the terminology ``K\"{a}hler'' for whenever the differential forms
have a Clifford algebra structure, regardless of their valuedness. In other
words, we mean ``K\"{a}hler'' whenever%
\begin{equation}
dx^{i}dx^{j}+dx^{j}dx^{i}=2\delta ^{ij}\text{ or }2\eta ^{ij},
\end{equation}%
is satisfied.

In section 4, we shall consider Clifford-valued K\"{a}hler differential
forms., i.e. the coexistence of the Clifford structures defined by (1) and a
tangent Clifford algebra of valuedness. The latter is defined by the relation%
\begin{equation}
\mathbf{a}_{i}\mathbf{a}_{j}\text{ }\mathbf{+}\text{ }\mathbf{a}_{j}\mathbf{a%
}_{i}=2\delta _{ij}\text{ or }2\eta _{ij},
\end{equation}%
where $\eta ^{ij}$ is constituted by ones and minus ones, depending on the
signature of the metric of the space considered. An appropriate term to
refer to its elements would be Clifford-valued clifforms. They belong to the
tensor product of the two algebras respectively defined by equations (1) and
(2). We shall not, however, use the term clifform, since we shall not be
interested in differential forms which are not so.

In section 5, we study the idempotents for the aforementioned symmetries.
They live in particular in a commutative algebras contained in that tensor
product structure. This change of the game leads to considering products of
three of the idempotents for one-parameter symmetries. This greater
availability of products is the prelude to a far greater interplay of
symmetries.

\section{Kaehler's treatment of translation and rotation symmetry}

We have brought into this section pertinent results of K\"{a}hler's calculus
(KC) of scalar-valued differential forms \cite{Kahler60}, \cite{Kahler61} %
\cite{Kahler62}. The last of these publications is a more comprehensive
treatment of the first one. Nevertheless, some important results are
exclusive to just one of those publications, like his treatment of Lie
differentiation \cite{Kahler60}, his treatment of symmetry in solutions of
differential systems \cite{Kahler61}, and his treatment of total angular
momentum \cite{Kahler62} (in the sense of linearly combining into one
expression all three components of angular momentum). All of these three
specific subjects will be reproduced in this section.

\subsection{Kaehler's monary and binary solutions with symmetry}

Monary idempotents are those belonging to a one-parameter symmetry. By the
term binary idempotents, we mean the product of two monary idempotents. A K%
\"{a}hler binary idempotent is constituted by the product of the idempotents
for third component of angular momentum and for time translation symmetry.
With a view to future physical applications, those idempotents will also be
our starting point for the construction of ternary idempotents, i.e.
constituted by the product of three monary idempotents. In the all too
important dimension four, it will not make much of a difference in which
order we approach the three factors. This will become obvious.

There are two main types of symmetry idempotents, translational and
rotational. Of special interest is the interplay between the two types and
also the linear coexistence of different one-parameter rotational
idempotents within a given equation. Such a coexistence happens in K\"{a}%
hler's work, where total angular momentum, though linear in its components
(in a generalized sense of the term) is not a vector operator.

Except for notational changes, K\"{a}hler defines the idempotents 
\begin{equation}
\epsilon ^{\pm }\equiv \frac{1}{2}(1\mp i\,dt),\;\ \ \ \ \ \ \;I_{xy}^{\pm
}\equiv \frac{1}{2}(1\pm i\,dx\,dy),
\end{equation}%
with $(dx^{i})^{2}=1=-dt^{2}$ and $c=1$. The inversion of the signs in the
definition of $\epsilon ^{\pm }$ has to do with his use of signature. The
unit imaginary is there for the purpose of making idempotents. The squares
of $\,idt$ and $i\,dx\,dy$ are $+1$, as required. The squares of $\,dt$ and $%
\,dx\,dy$ are not.

The $\epsilon ^{\pm }$ commute with the $I_{xy}^{\pm }$. We also have 
\begin{equation}
\epsilon ^{+}+\epsilon ^{-}=1=I_{xy}^{+}+I_{xy}^{-},\;\ \ \ \ \ \;\epsilon
^{\pm }\epsilon ^{\mp }=I_{xy}^{\pm }I_{xy}^{\mp }=0.
\end{equation}%
From these equations follows that the statement 
\begin{equation}
u=u^{+}\text{ }\epsilon ^{+}\text{ }+\text{ }u^{-}\text{ }\epsilon ^{-}
\end{equation}%
defines $dt$-independent \textit{spatial differentials} $u^{\pm }.$ Given a
spatial differential, $u$, the statement 
\begin{equation}
u=\text{ }^{_{+}}u\text{ }I_{xy}^{+}+\text{ }^{_{-}}u\text{ }I_{xy}^{-}
\end{equation}%
similarly defines so called \textit{meridian differentials} $^{\pm }u$, i.e.
that they do not depend on $d\phi $ when written in terms of the basis ($%
d\rho ,d\phi $, $dz$). The coefficients in (5)-(6) may still depend on all
coordinates.

Combining these definitions and observations, we have 
\begin{equation}
1\equiv \epsilon ^{+}I_{xy}^{+}+\epsilon ^{+}I_{xy}^{-}+\epsilon
^{-}I_{xy}^{+}+\epsilon ^{-}I_{xy}^{-},
\end{equation}%
and, therefore,%
\begin{equation}
u\equiv u\epsilon ^{+}I_{xy}^{+}+u\epsilon ^{+}I_{xy}^{-}+u\epsilon
^{-}I_{xy}^{+}+u\epsilon ^{-}I_{xy}^{-}.
\end{equation}%
The set of the four idempotents $\epsilon ^{\pm }I_{xy}^{\ast }$ spans the
same 4-dimensional module over the complex numbers as the set $%
(1,dt,dxdy,dtdxdy)$ does. The asterisk means that we may choose for $I_{xy}$
the superscripts $+$ and $-$ independently of the specific superscript in $%
\epsilon .$

\subsection{Of Kaehler's idempotents and commutativity}

Idempotents like $\epsilon ^{\pm }$ require only a $1-$dimensional space. We
also need just a one-dimensional space for space translation idempotents. $%
I_{xy}^{\pm }$ requires at least but not necessarily more than a $2-$%
dimensional space, where we also have $(1/2)(1\pm i\,dx)$ and $(1/2)(1\pm
i\,dy).$ But the last two ones, when not mutually annulling, do not commute,
nor do they commute with the rotational ones (Commutativity will be key in
the considerations that follow). In two dimensions, we could also have a ($%
t,x$) space, the terminology speaking of the signature of the metric. We
will not deal, however, with hyperbolic rotation symmetry, as it will loose
relevance when we deal with Lorentzian spaces of dimensions three and four.

The idempotents in (3) would fit in $3-$dimensional and $4-$dimensional
spaces, ($t,x,y$) and ($t,x,y,z$). Why did K\"{a}hler not consider
translations in the $x$, $y$ and $z$ dimensions in his study of symmetry in
spacetime? Once he had to consider $I_{xy}^{\pm }$ in order to deal with
rotations for the problem of the fine structure of the hydrogen atom (and
also of time translation because of rest mass), there is no room for $%
(1/2)(1\pm i\,dx^{i})$ since they do not commute with $\epsilon ^{\pm }.$
And there is no room for $I_{yz}^{\pm }$ and $I_{zx}^{\pm }$ (once $%
I_{xy}^{\pm }$ has been chosen), since they do not commute with $I_{xy}^{\pm
}.$

We proceed to see what role commutativity and the mutually annulling of
idempotents exhibited in the last equations (4) play. The impact of
commutativity is that it allows for the mutually annulling of idempotents
such as, say, $I_{xy}^{+}\epsilon ^{-}$ and $I_{xy}^{-}\epsilon ^{-}.$
Suppose we multiply (7) by $I_{xy}^{-}\epsilon ^{-}$ on the right. We get $%
I_{xy}^{+}\epsilon ^{-}I_{xy}^{-}\epsilon ^{-}=I_{xy}^{+}I_{xy}^{-}\epsilon
^{-}\epsilon ^{-}=0.$ If $I_{xy}^{+}$ and $\epsilon ^{-}$ did not commute,
we could not have obtained this result. The impact of the vanishing of three
of the four terms annulling when multiplying on the right of (7) by each of
the four $\epsilon ^{\pm }I_{xy}^{\ast }$ can be ascertained by trying to
make the same considerations with 
\begin{equation}
1\equiv
P_{x}^{+}I_{xy}^{+}+P_{x}^{+}I_{xy}^{-}+P_{x}^{-}I_{xy}^{+}+P_{x}^{-}I_{xy}^{-},
\label{9}
\end{equation}%
and, therefore,%
\begin{equation}
u\equiv
uP_{x}^{+}I_{xy}^{+}+uP_{x}^{+}I_{xy}^{-}+uP_{x}^{-}I_{xy}^{+}+uP_{x}^{-}I_{xy}^{-},
\end{equation}%
where $P_{x}^{\pm }\equiv \frac{1}{2}(1\pm i\,dx)$. These are valid
equations but the terms on the right of (9) and (10) do not have the same
properties as the terms on the right of (7) and (8), which are disjoint,
read orthogonal.

The $\epsilon ^{\pm }I_{xy}^{\ast }$ are so call constant differentials \cite%
{Kahler60}, \cite{Kahler62}, meaning $\partial (\epsilon ^{\pm }I_{xy}^{\ast
})=0$, where $\partial $ is the sum of the interior (read also divergence
and co-derivative) and exterior derivatives. As a consequence of one of the
properties of constant differentials, we have%
\begin{equation}
\partial (u\epsilon ^{\pm }I_{xy}^{\ast })=(\partial u)\epsilon ^{\pm
}I_{xy}^{\ast }.
\end{equation}%
Hence all four terms on the right hand side of (8) are solutions of the same
equation of the type $\partial u=au$ if $u$ is. The $u\epsilon ^{\pm
}I_{xy}^{\ast }$ are mutually orthogonal solutions of the same equation. We
are then justified in writing%
\begin{equation}
u=\text{ }^{+}u^{+}\text{ }I_{xy}^{+}\epsilon ^{+}\text{ }+\text{ }^{+}u^{-}%
\text{ }I_{xy}^{+}\epsilon ^{-}\text{ }+\text{ }^{-}u^{+}\text{ }%
I_{xy}^{-}\epsilon ^{+}\text{ }+\text{ }^{-}u^{-}\text{ }I_{xy}^{-}\epsilon
^{-},
\end{equation}%
with uniquely defined meridian differentials $^{\pm }u^{\ast }.$

One is thus led to seek spinorial solutions that take the form of meridian
differentials times $I_{xy}^{\pm }\epsilon ^{\ast }$. In a exterior system,
the dependence on $t$ and $\phi $ of proper functions for time translations
and rotational symmetry takes place through the standard phase factors. It
is thus \textit{unavoidable} (English translation of the term used by K\"{a}%
hler to make the present point) to consider spinorial solutions of the form 
\begin{equation}
u=pe^{im\phi -iEt}\epsilon ^{\pm }I_{xy}^{\ast }
\end{equation}%
where $p$ is a \textit{strict meridian differential, }i.e. meridian
differential whose coefficients depend on $\rho $ and $z$, but not on $\phi $
and $t$.

With a view to future developments, it is worth pointing out at the nature
of the different factors in the last expression. The $e^{im\phi -iEt}$ and $%
\epsilon ^{\pm }I_{xy}^{\ast }$ are structural. On the other hand, $p$ will
contain more fleeting info, not further restricted by those symmetry
considerations but by the specifics of particular problems. Between $%
e^{im\phi -iEt}$ on the one hand and $\epsilon ^{\pm }I_{xy}^{\ast }$ on the
other hand, the idempotents contain the most rigid factor in (13); in the
phase factor, one still has the freedom of having different values for $m$
and $E.$ There is not a comparable freedom in $\epsilon ^{\pm }I_{xy}^{\ast
} $.

Let us finally deal with the relation of $\epsilon ^{\pm }$ to charge, in
addition to energy. In 1962, K\"{a}hler derived a conservation law%
\begin{equation}
d(u,\eta \bar{u})_{1}=0,  \label{14}
\end{equation}%
from his equation that replaces Dirac's \cite{Kahler62}. It is quadratic, in
the sense that it is bilinear in $u$ and its complex conjugate, $\bar{u}$.
The operator $\eta $ reverses the sign of differential forms of odd grade,
and $(\_,\_)_{1}$ is a product of differential forms of grade $n-1$, where $%
n $ is the dimensionality of the space. We need not enter more details for
present purposes. It is well known that equations like this represent
conservation laws. What is conserved in this case?

K\"{a}hler was able to develop (14) to the point where it acquired the form%
\begin{equation}
\frac{\partial \rho }{\partial t}+j+\frac{\partial \rho ^{\prime }}{\partial
t}+j^{\prime }=0
\end{equation}%
(Needless to say that $\rho $ and $\rho ^{\prime }$ stand here for charge
densities, and $j$ and $j^{\prime }$ stand for physical currents). Primed
and unprimed quantities correspond to the two terms on the right of (5),
which associates $\epsilon ^{\mp }$ with the two signs of charge, not of
energy, which is the same in both cases. The problem of associating
particles with antiparticles with opposite sign of energy does not exist in K%
\"{a}hler's work.

\subsection{Kaehler's approach to Lie differentiation}

We have dealt with the form of solutions with the stated symmetries. We now
deal with the related issue of operators for the same symmetries.

In KC, spin emerges through Lie differentiation as an integral part of
angular momentum, not as a late attachment to the orbital one. One has to
ask oneself, like K\"{a}hler did, what is the action of an operator%
\begin{equation}
A\equiv \xi _{i}(x^{1},...,x^{n})\frac{\partial }{\partial x^{i}}  \label{16}
\end{equation}%
on differential forms. It is not what would appear to be at first sight.
Suffice to consider the case of definite grade:%
\begin{equation}
u\equiv \frac{1}{p!}a_{i_{1}...i_{p}}dx^{i_{1}}\wedge ...\wedge dx^{i_{p}}.
\end{equation}%
The issue thus lies in extending the action of operator (16) from the ring
of differential $0-$forms to the ring of differential forms. K\"{a}hler (%
\cite{Kahler62}, section 16) rightly credits Cartan with such an extension %
\cite{Cartan22}. But the latter's treatment is too cryptic, as he was not
interested in providing an expression for the action of the operator. He
instead directly derived with a very brief argument the theorem that is
sometimes referred to as the golden theorem. It amounts to a different
version of that action. Because of its clarity, we now reproduce K\"{a}%
hler's derivation of the formula for $Au$ that he published in \cite%
{Kahler60}.

The right hand side of (16) could be viewed as the partial derivative with
respect to some coordinate in suitable coordinate systems ($y^{i}$). In
order to find them, K\"{a}hler considers the general solution of the system%
\begin{equation}
\frac{dx^{i}}{dy^{n}}=\xi ^{i}(x^{1},...,x^{n}),
\end{equation}%
which one may express as%
\begin{equation}
x^{i}=x^{i}(y^{1},...,y^{n}),
\end{equation}%
where the $y^{1},...,y^{n-1}$ are first integrals of (18)$.$ The last of the
first integrals is additive to the variable $y^{n}$, with which we may thus
identify it. In the neighborhood of a point, (19) can be viewed as a
coordinate transformation. The pull-back of (17) to the $y$ coordinate
system is%
\begin{equation}
u\equiv \frac{1}{p!}a_{i_{1}...i_{p}}\frac{\partial x^{i_{1}}}{\partial
y^{k_{1}}}\text{...}\frac{\partial x^{i_{p}}}{\partial y^{k_{p}}}%
dy^{k_{1}}\wedge ...\wedge dy^{k_{p}}.
\end{equation}

We now compute $\frac{\partial u}{\partial y^{n}}$ using that%
\begin{equation}
\frac{\partial }{\partial y^{n}}=\frac{\partial x^{i}}{\partial y^{n}}\frac{%
\partial }{\partial x^{i}}=\xi ^{i}\frac{\partial }{\partial x^{i}}=A.
\label{21}
\end{equation}%
Notice that, in (18), we had the $\xi ^{i}$ as total derivatives since there
was just one independent variable. The same variable is viewed as the last
one in a different coordinate system, hence the partial differentiation
notation. The derivatives continue to be the same. The result of the
computation is%
\begin{equation}
Au=\xi ^{i}\frac{\partial u}{\partial x^{i}}+d\xi ^{i}\wedge e_{i}u.
\end{equation}%
K\"{a}hler writes $d\xi ^{i}$ as $d(\xi ^{i})$ since he reserves the
notation $d\xi ^{i}$ for covariant derivatives (See \cite{Kahler60}, end of
his section 17, and second equation in section 34). The non-intuitive
implication%
\begin{equation}
Au\neq \frac{1}{p!}\xi _{i}\frac{\partial a_{i_{1}...i_{p}}}{\partial x^{i}}%
dx^{i_{1}}\wedge ...\wedge dx^{i_{p}}
\end{equation}%
is due to the fact that what remains constant in partial differentiations
changes from term to term in (16). In order to avoid inadvertent errors, one
should thus view $\partial /\partial \phi ^{i}$ as the starting point for
rotational symmetry with preference over 
\begin{equation}
J_{i}=x^{j}\frac{\partial }{\partial x^{k}}-x^{k}\frac{\partial }{\partial
x^{j}},  \label{24}
\end{equation}%
Here $\phi ^{i}$ represents the angle of rotation around three orthonormal
axes and ($i,j,k$) constitute cyclic permutations of ($1,2,3$). Clearly each 
$\phi $ plays the role of $y^{n}$ in the previous argument. Coordinates $%
y^{1},...,y^{n-1}$ would then be the other coordinates in systems such as
those of spherical and cylindrical coordinates, but different systems for
different index $i$ in $\phi ^{i}$.

\subsection{Angular momentum in Kaehler's theory}

In terms of Cartesian coordinates, the action of $\partial /\partial \phi
^{i}$ on differential forms is given by%
\begin{equation}
J_{i}u=x^{j}\frac{\partial u}{\partial x^{k}}-x^{k}\frac{\partial u}{%
\partial x^{j}}+\frac{1}{2}w_{i}u-\frac{1}{2}uw_{i},
\end{equation}%
where%
\begin{equation}
w_{i}=dx^{j}dx^{k}\equiv dx^{jk}.
\end{equation}%
\cite{Kahler60}, \cite{Kahler62}. We shall refer to $\frac{1}{2}w_{i}u-\frac{%
1}{2}uw_{i}$ as the spin term(s), though chirality terms would also be
appropriate. These two terms are intrinsically associated, like (we have
seen) energy and charge also are.

Notice that, if $w_{i}$ and $u$ commute, the last two terms in (25) add up
to zero. Such $u$'s are proper values of $\frac{1}{2}w_{i}\_\_-\frac{1}{2}%
\_\_w_{i}$ with proper value zero. The same two terms add up to $w_{i}u$ if
they anticommute. Clearly, $w_{k}$ commutes with $dx^{k}$ and $%
(1/2)(1+w_{k}) $. It anticommutes with $dx^{j}$, $dx^{k}$, $w_{j}$ and $%
w_{k} $, and it neither commutes nor anticommutes with $(1/2)(1+w_{i})$ and $%
(1/2)(1+w_{j})$.

Consider now total angular momentum, not just components thereof. K\"{a}hler
defines an operator $K$ by%
\begin{equation}
(K+1)u\equiv \sum_{i=1}^{i=3}J_{i}u\vee w_{i},
\end{equation}%
acting on scalar-valued differential forms $u.$ The reason for K\"{a}hler's
defining $(K+1)u$ where he could have simply written $Ku$ has the unstated
purpose of later ending with the same terminology as in the standard
literature. He showed that%
\begin{equation}
(K+1)^{2}u=\sum_{i=1}^{i=3}J_{i}^{2}u+(K+1)u.
\end{equation}%
If we take the last term to the left side, we get the equation%
\begin{equation}
(K+1)Ku=\sum_{i=1}^{i=3}J_{i}^{2}u.  \label{29}
\end{equation}%
It is a relation about operators which has a parallel in a proper values
equations in the standard theory of angular momentum. Equation (27)
implicitly speaks of $J_{i}$ as its components, though not in the standard
sense of components in a vector space or in a module. The operator 
\begin{equation}
K+1\equiv \sum_{i=1}^{i=3}J_{i}\_\vee w_{i}
\end{equation}%
is differential $2-$form valued, in the sense that it is the contraction of
the $J_{i}$ with the elements of a basis of differential $2-$forms.

The issue of the grade of the operator $K+1$ is a subtle one. We write its
components as 
\begin{equation}
J_{i}\equiv =x^{j}\frac{\partial }{\partial x^{k}}-x^{k}\frac{\partial }{%
\partial x^{j}}+\frac{1}{2}w_{i}\vee \_\_-\frac{1}{2}\_\_\vee w_{i},
\label{31}
\end{equation}%
where we have made explicit the symbol for Clifford product. The first two
terms on the right of (31) are of grade zero; the last two appear to be of
grade two. This is, however, deceiving, as we now explain.

Start by considering that the operator definition (31) follows (25), not the
other way around. Given arbitrary differential forms $u$ and $v,$ each of
homogeneous grade, the grade of $u\vee v-v\vee u$ is two units less than the
sum of the respective grades (see \cite{Kahler62}, p. 478). By virtue of
this, the grade of $J_{i}u$ and $(K+1)u$ is the same as the grade of $u$ (if
of definite grade $r$) in spite of the fact that each of $\frac{1}{2}w_{i}u$
and $\frac{1}{2}uw_{i}$ may consist of terms of grades $r+2,$ $r$ and $r-2$.

Notice that, in this treatment of angular momentum, there has not been any
need for the unit imaginary. At a certain point, however, it becomes
necessary in K\"{a}hler's treatment in order to match operators with
idempotents and thus with phase factors.

For completeness purposes, we call attention to his remark that the
components of angular momentum should be $\frac{h}{2\pi i}J_{i}$, ``and not
as is usually the case%
\begin{equation}
\frac{h}{2\pi i}\left( J_{1}+\frac{1}{2}\text{ \ }\vee w_{1}\right) ,\text{
\ \ \ }\frac{h}{2\pi i}\left( J_{2}+\frac{1}{2}\text{ \ }\vee w_{2}\right) ,%
\text{ \ \ \ }\frac{h}{2\pi i}\left( J_{3}+\frac{1}{2}\text{ \ }\vee
w_{3}\right)
\end{equation}%
\cite{Kahler62} (section 27). (Compare with \cite{Gilmore}, chapter 7,
Exercise 10).

\subsection{Action of angular momentum on rotation idempotents}

The action of $J_{l}$ on $w_{k}$ and on the idempotents $\frac{1}{2}(1\pm
iw_{k})$ will be the same as the action of $\frac{1}{2}w_{l}\_\_-\frac{1}{2}%
\_\_w_{l}.$ Obviously%
\begin{equation}
J_{k}I_{ij}^{\pm }=\frac{1}{2}w_{k}\frac{1}{2}(1\pm iw_{k})-\frac{1}{2}\frac{%
1}{2}(1\pm iw_{k})w_{k}=0.
\end{equation}%
The idempotents $I_{ij}^{\pm }$ are proper functions of $J_{k}$ and $\frac{1%
}{2}(w_{k}\_\_-\_\_w_{k})$ with proper value zero. Nothing impedes that spin
still be included in the factor $e^{im\phi }$ in (13), through the first two
terms in (25). The proper value zero becomes a pair of opposite proper
values for solutions which have $I_{ij}^{\pm }$ as factor and anticommute
with the operator $J_{k}$. The difference $\frac{1}{2}w_{k}u-\frac{1}{2}%
uw_{k}$ becomes $-uw_{k}.$ Then $w_{k}$ is absorbed into $(1\pm iw_{k}).$
The sign of the proper value will depend on the sign of the action of $w_{k}$
on $I_{ij}^{\pm }.$

For the other combinations of indices, we first get%
\begin{equation}
J_{i}w_{k}=w_{j},\text{ \ \ \ \ \ \ \ \ \ \ }J_{j}w_{k}=-w_{i}.  \label{34}
\end{equation}%
Then, clearly,%
\begin{equation}
J_{i}I_{ij}^{\pm }=\pm \frac{1}{2}iw_{j},\text{ \ \ \ \ }J_{i}I_{ij}^{\pm
}=\mp \frac{1}{2}iw_{i}.
\end{equation}%
The $I_{ij}^{\pm }$ are not proper ``functions'' of $J_{i}$, $J_{j}$, $\frac{%
1}{2}(w_{i}\_\_-\_\_w_{i})$ and $\frac{1}{2}(w_{j}\_\_-\_\_w_{j}).$

We use (34) and (33) to obtain%
\begin{equation}
\sum_{i=1}^{i=3}J_{l}w_{1}\vee
w_{l}=0+w_{3}w_{2}-w_{2}w_{3}=2w_{3}w_{2}=2w_{1}.
\end{equation}%
Thus, in general,%
\begin{equation}
(K+1)w_{k}=2w_{k}.
\end{equation}%
Hence, $w_{k}$ and $iw_{k}$ are proper functions of $K$ with proper value $%
+1.$ Finally,%
\begin{equation}
(K+1)\frac{1}{2}(1\pm iw_{k})=\pm (K+1)iw_{k}=\pm iw_{k}.
\end{equation}%
The rotational idempotents are not per se proper functions of $K$, but can,
of course, make part of expressions that are. Let us solve for $K\frac{1}{2}%
(1+iw_{k})$. We have, using (38),%
\begin{equation}
K\frac{1}{2}(1\pm iw_{i})=\pm iw_{i}-\frac{1}{2}(1\pm iw_{i})=-\frac{1}{2}%
\pm \frac{1}{2}iw_{i}=-\frac{1}{2}(1\mp iw_{i}).
\end{equation}%
It thus follows that $-K$ moves $\frac{1}{2}(1\pm iw_{i})$ down and up:%
\begin{equation}
-K\frac{1}{2}(1\pm iw_{i})=\frac{1}{2}(1\mp iw_{i}).
\end{equation}

Notice that the unit imaginary has come into play through the idempotents.
It had not done so through the angular momentum operators on their own.

\section{Symmetry Related Transformations and the\protect\linebreak\ Unit
Imaginary}

In this section, we shall treat rotations through Clifford algebra, which is
the canonical algebra of Euclidean and pseudo-Euclidean spaces. As a
by-product, we associate rotational symmetries with gauge transformations.
This association can then be extended to other symmetries. All this is
relevant for our later replacement of the unit imaginary in the next
section. It usurps the role of geometric quantities of square $-1.$

\subsection{Rotation of vectors}

The contents of this subsection is well known by every practitioner of
Clifford algebra and can be found in one of our posted papers in arXiv \cite%
{V50}. We repeat it here for completeness purposes. In the plane, let $%
\mathbf{t}$ be a unit vector, and let $\mathbf{u}^{\prime }$ be the
reflection of $\mathbf{u}$ with respect to the direction $\mathbf{t}$. We
clearly have $\mathbf{u}^{\prime }\mathbf{t=tu}$ because ($\mathbf{u}%
^{\prime },\mathbf{t}$) is indistinguishable from ($\mathbf{t,u}$) (picture
this situation with $\mathbf{u}$ of unit magnitude). Since $\mathbf{t}^{-1}=%
\mathbf{t}$, we have%
\begin{equation}
\mathbf{u}^{\prime }\mathbf{=tut}^{-1}\text{ }\mathbf{=tut.}  \label{41}
\end{equation}%
Consider next the reflection of $\mathbf{u}$ with respect to a plane $\pi $
with normal unit vector $\mathbf{n.}$ Let $\mathbf{t}$ now be a unit vector
in the direction of the projection of $\mathbf{u}$ on $\pi .$ The three
vectors ($\mathbf{n}$, $\mathbf{t}$ and $\mathbf{u}$) are on the same plane, 
$\mathbf{n}$ and $\mathbf{t}$ being perpendicular. By a known elementary
theorem, the reflections of $\mathbf{u}$ with respect to two perpendicular
vectors are opposite:%
\begin{equation}
\mathbf{nun=-tut.}
\end{equation}

In 3-D, a rotation of $\mathbf{u}$ by an angle $\phi $ around an axis in the
direction of the unit vector $\mathbf{n}_{3}$ yields the vector $\mathbf{u}%
^{\prime \prime }$ resulting from successive reflections on two planes
making an angle $\phi /2$, intersecting along $\mathbf{N}$ and with normal
unit vectors $\mathbf{n}_{1}$ and $\mathbf{n}_{2}.$ We thus have%
\begin{equation}
\mathbf{u}^{\prime \prime }\text{ }\mathbf{=n}_{2}\mathbf{n}_{1}\mathbf{un}%
_{1}\mathbf{n}_{2}=(\mathbf{n}_{1}\mathbf{n}_{2})^{-1}\mathbf{u(n}_{1}%
\mathbf{n}_{2})\mathbf{.}
\end{equation}%
We write the Clifford product $\mathbf{n}_{1}\mathbf{n}_{2}$ in terms of its
scalar and exterior components:%
\begin{equation}
\mathbf{n}_{1}\mathbf{n}_{2}\mathbf{=}\cos \frac{\phi }{2}+\lambda \mathbf{a}%
_{1}\mathbf{a}_{2},
\end{equation}%
where $\mathbf{a}_{1}$ and $\mathbf{a}_{2}$ are any two orthonormal vectors
in the plane determined by $\mathbf{n}_{1}$ and $\mathbf{n}_{2}.$ Clearly
then%
\begin{equation}
\mathbf{n}_{1}\mathbf{n}_{2}\text{ }\mathbf{=}\cos \frac{\phi }{2}+\sin 
\frac{\phi }{2}\mathbf{a}_{1}\mathbf{a}_{2}=e^{\frac{\phi }{2}\mathbf{a}_{1}%
\mathbf{a}_{2}},
\end{equation}%
since $(\mathbf{a}_{1}\mathbf{a}_{2})^{2}=-1.$ It thus follows that%
\begin{equation}
\mathbf{u}^{\prime \prime }\text{ }\mathbf{=}\text{ }e^{-\frac{\phi }{2}%
\mathbf{a}_{1}\mathbf{a}_{2}}\text{ }\mathbf{u}\text{ }e^{\frac{\phi }{2}%
\mathbf{a}_{1}\mathbf{a}_{2}}.
\end{equation}

\subsection{(Mainly) rotational and related gauge symmetry}

We now extend the action of rotations to elements of the tangent Clifford
algebra. We apply it to monomials; the extension\ to polynomials is trivial.
Define%
\begin{equation}
\mathbf{R\equiv }\text{ }e^{\frac{\phi }{2}\mathbf{a}_{i}\mathbf{a}_{j}}.
\label{47}
\end{equation}%
Under a rotation around $\mathbf{a}_{k}$, the Clifford product of vectors $%
\mathbf{A\equiv ab...k}$ becomes%
\begin{equation}
\mathbf{A}^{\prime }\mathbf{\equiv (R}^{-1}\mathbf{aR)\mathbf{(R}^{-1}%
\mathbf{bR)}...\mathbf{(R}^{-1}\mathbf{kR)=\mathbf{R}^{-1}\mathbf{a\mathbf{b}%
}}}...\mathbf{\mathbf{\mathbf{\mathbf{kR=}R}^{-1}\mathbf{AR}.}}
\end{equation}%
Spinors are members of left ideals. Rotational spinors in the tangent
algebra are the elements of ideals%
\begin{equation}
\mathbf{A}\frac{1}{2}(1\pm i\mathbf{a}_{i}\mathbf{a}_{j}),
\end{equation}%
where $\mathbf{A}$ represents all elements in the algebra. The unit
imaginary is there to build the idempotents since ($\mathbf{a}_{i}\mathbf{a}%
_{j}$)$^{2}$ equals $-1$, and not the required $+1$. Hence something is not
quite right with Clifford algebra, as one should not need units imaginary to
deal with idempotents in real vector spaces. It looks as if there must be
``something Clifford'' which is more fundamental than the Clifford algebra
itself. We have the solution in the tensor product of two specific, related
Clifford algebras. But that is not the only reason why we shall have
something more fundamental \cite{Kahler60}. More on this in the next section.

Since (49) belongs to the algebra, its rotation in the plane determined by $%
\mathbf{a}_{i}$ and $\mathbf{a}_{j}$ is given by 
\begin{equation}
\mathbf{R}^{-1}\mathbf{A}\frac{1}{2}(1\pm i\mathbf{a}_{i}\mathbf{a}_{j})%
\mathbf{R}=\mathbf{R}^{-1}\mathbf{AR}\frac{1}{2}(1\pm i\mathbf{a}_{i}\mathbf{%
a}_{j})=\mathbf{A}^{\prime }\frac{1}{2}(1\pm i\mathbf{a}_{i}\mathbf{a}_{j}),
\label{50}
\end{equation}%
which means that the ideal transforms into itself. But general members of
this ideal transform to $\mathbf{R}^{-1}\mathbf{A}\frac{1}{2}(1\pm i\mathbf{a%
}_{i}\mathbf{a}_{j})\mathbf{R}$, not to%
\begin{equation}
\mathbf{R}^{-1}\mathbf{A}\frac{1}{2}(1\pm i\mathbf{a}_{i}\mathbf{a}_{j}).
\end{equation}%
Transformation (51) ---equivalently ($\mathbf{\mathbf{R}^{-1}A\mathbf{R}}$, $%
\mathbf{\mathbf{R}^{-1}\psi }$)--- is the gauge transformation associated
with but not identical to a rotation, even if $\mathbf{R}$ is given by (47).

By virtue of their structure, gauge transformations apply to equations of
the type%
\begin{equation}
(\mathbf{\partial -A)\psi =0,}
\end{equation}%
where, suffice to say, $\mathbf{A}$ and $\mathbf{\psi }$ are respective
input and output and where $\mathbf{R}$ might not be as for rotations, (47).
The exponent in $\mathbf{R}$ could be, for instance, a scalar times a vector
of square $+1$, or $i$ times a vector of square $-1$. This is so by virtue
of the fact that 
\begin{equation}
\lbrack \mathbf{R}^{-1}(\mathbf{\partial -a)\mathbf{R]}}\text{ }\mathbf{%
\mathbf{R}^{-1}\psi =0}
\end{equation}%
follows from (52) for any operator which has an inverse. Of course, the
physical relevance of $\mathbf{R}$ will depend on whether $\mathbf{R}^{-1}(%
\mathbf{\partial -a)\mathbf{R}}$ is an invariant.

Let $\lambda $ be an arbitrary scalar and let $\mathbf{R}$ be of the type $%
e^{\lambda \mathbf{\Sigma }}$ with $\mathbf{\Sigma }^{2}=-1.$ The
idempotents 
\begin{equation}
\frac{1}{2}(1\pm i\mathbf{\Sigma })
\end{equation}%
define ideals%
\begin{equation}
\mathbf{A}\frac{1}{2}(1\pm i\mathbf{\Sigma })  \label{55}
\end{equation}%
associated with the gauge transformation%
\begin{equation}
(\mathbf{\partial -a)\rightarrow R}^{-1}(\mathbf{\partial -a)\mathbf{R}}%
\text{, \ \ \ \ \ \ \ \ }\mathbf{\psi \rightarrow R}^{-1}\mathbf{\psi .}
\end{equation}%
If the two idempotents (54) are associated with solutions of an equation
possessing some one-parameter group of classical geometric symmetry, the
gauge transformation (56) is associated with that symmetry and we may say
that it is a classical gauge symmetry, as opposed to gauge symmetries not
directly related to the tangent bundle.

We have seen gauge transformations as invariance properties of equations,
(56). If $\mathbf{R}$ is as in (50),the second of these can then be seen as
a gauge transformation of the associated ideal, in turn associated with a
rotation. Hence the two groups in transformations (56) (rotation groups and
spin groups) respond to the same symmetry. Hence, it is correct but not
illuminating to speak of two groups. They reflect different applications of
the same symmetry.

\subsection{Of units imaginary in translational and related gauge symmetry}

Assume that K\"{a}hler had been concerned with gauge symmetry under space
translation. Instead of the $\epsilon ^{\pm }$ given in (3), he would have
had idempotents like $\frac{1}{2}(1\pm dz).$ Then $\mathbf{R}$ would be $%
\mathbf{R\equiv }$ $e^{-\lambda \mathbf{a}_{3}}$. It is not relevant here
what $\mathbf{\mathbf{R}^{-1}...\mathbf{R}}$ would mean. We have said gauge
symmetry under space translation, not space translation symmetry. This
difference, which must be clear from the contents of the previous
subsection, may be relevant in quantum mechanics. It may be the beginning of
a hint of why one has problems with position operators in standard quantum
mechanics.

Similarly, we would not have a phase factor for hyperbolic rotation
symmetry. The idempotent would be like $\frac{1}{2}(1\pm dtdz)$, and the
square of $dtdz$ is plus one.

We shall see in the next section that the right choice of algebra makes the
unit imaginary unnecessary for building idempotents or for gauge
transformations. We may thus advance that the presence/absence of the
standard unit imaginary is a function of the algebra. From this point on in
this paper, we shall deal with tangent Clifford algebra valued differential
forms. Hence, the gauge symmetries that we shall explicitly or implicitly
consider will be of the geometric type, meaning that they pertain to the
tangent bundle of vectors, of Clifford elements, of frames, of spinors, etc.
These are called tangent, Clifford, frame and spinor bundles in the
literature. Here, the term tangent bundle will refer to all of them, in
contradistinction to auxiliary bundles.

The natural question then is whether a non-geometric gauge transformation
actually is a geometric one that simply has gone unnoticed by an unfortunate
choice of algebra. But that is a physical issue, which will thus not be
considered in this paper. Let us summarize by saying that there is no need
for the standard unit imaginary when one deals with the algebras that shall
occupy us in the next\ section. The issue of the unit imaginary is far more
subtle than usually acknowledged. Here is a summary.

\section{Clifford-Valued Clifforms}

An algebraically minded mathematician who were exposed for the first time to
the expression $dx\mathbf{i+}dy\mathbf{\mathbf{j+}}dz\mathbf{k}$ might see
in there the tensor product of two modules, and/or the tensor product of the
two exterior or Clifford algebras spanned by those modules. On the other
hand, most of us have been introduced to the wrong concept of ($dx\mathbf{,}%
dy,dz$) so that $dx\mathbf{i+}dy\mathbf{\mathbf{j+}}dz\mathbf{k}$ is a
vector with weird components. It is rather a vector-valued differential
form, meaning a vector-valued function of curves. It was already obvious to 
\'{E}. Cartan that one should be dealing with tensor products of structures;
he referred to the curvature of Euclidean connections as bivector-valued
differential 2-forms \cite{Cartan23} (See specially his sections on integral
invariants in the chapter on Riemannian manifolds with metric connections of
that reference). Both of the structures in his product were exterior
algebras. We shall be interested in viewing $dx\mathbf{i+}dy\mathbf{\mathbf{%
j+}}dz\mathbf{k}$ as pertaining to the tensor product of two Clifford
algebras.

\subsection{Tensor product of Clifford algebras}

Exterior and Clifford algebras are quotient algebras of general tensor
algebras. Being algebras, they are modules or vector spaces in the first
place, of dimension $2^{n}$ if constructed upon modules or vector spaces of
dimension $n$. Consider the tensor product of the two Clifford algebras
respectively defined by the relations (1) and (2) in terms of orthonormal
bases. Simple bases in those algebras are ($1,dx$, $dy$, $dz$, $dydz$, $dzdx$%
, $dxdy$, ...) and ($1,$ $\mathbf{i}$, $\mathbf{j}$, $\mathbf{k}$, $\mathbf{%
jk}$, $\mathbf{ki}$, $\mathbf{ij}$, $\mathbf{jk}$, ...\textbf{) }%
respectively.\textbf{\ }The element $dx\mathbf{j}$ of the product structure
has components ($0,$ $\mathbf{j}$, $0$, $0$, $0$, $0...$) in one algebra and
($0,0,dx,0,0,0,...$) in the other. A basis in the product structure is
constituted by the $2^{n}\times 2^{n}$ elements resulting from product of
all the basis elements in one algebra by all the basis elements in the other
one.

Each Clifford algebra factor is isomorphic to the structure spanned by the
product of ($1,dx$, $dy$, $dz$, $dydz$, $dzdx$, $dxdy$, ...) with ($1,$ $0$, 
$0$, ...$0$). The other Clifford algebra factor is isomorphic to the product
of ($1,$ $0$, $0$, ...$0$) with ($1,$ $\mathbf{i}$, $\mathbf{j}$, $\mathbf{k}
$, $\mathbf{jk}$, $\mathbf{ki}$, $\mathbf{ij}$, $\mathbf{jk}$, ...\textbf{)}

In the new structures just considered, we have Clifford-valued clifforms,
the term clifform being used when the exterior product of differential forms
is replaced with their Clifford product. K\"{a}hler developed calculus for
tensor-valued differentials (clifforms, which he called tensors whose
coefficients are differential forms. And he\ considered physical
applications for scalar-valued differential forms, but not for those of more
general valuedness. In particular, he did not extend to them his study of
symmetry in solutions of exterior systems.

In more recent times, tensor products of Clifford algebras have been
considered by Oziewicz \cite{Oziewicz}. Interestingly, products of similar
type have emerged in very deep algebraic work by Helmstetter \cite{Helms}.
Unlike Oziewicz'es and our own work, such tensor products are not per se the
primary concepts. In view of the work in this paper, we submit that, ``in
matters Clifford'', there is something more fundamental than a single
Clifford algebra, whether based on equations such as (1) or such as (2).

\subsection{Commutative, mirror algebra}

The $2^{2n}$ structure considered in the previous subsection has a subset of
mirror-symmetric elements, which give rise to a very interesting structure.
We proceed to explain the concept via some simple examples. $dx^{1}\mathbf{a}%
_{1}$ is mirror symmetric and $dx^{1}\mathbf{a}_{2}$ is not, since the
superscript of one factor is not equal to the subscript of the other factor.
But $dxdy\mathbf{ij}$ is mirror symmetric and $dxdy\mathbf{\mathbf{ik}}$ is
not. Certainly, $d\mathbf{r\equiv }$ $dx\mathbf{i+}dy\mathbf{j+}dz\mathbf{k}$
is a non-monomial example of mirror symmetric element. The inhomogeneous
grade element $dx\mathbf{i+}dxdy\mathbf{\mathbf{ij}}$ also is. As a member
of the structure spanned by ($1,$ $\mathbf{i}$, $\mathbf{j}$, $\mathbf{k}$, $%
\mathbf{jk}$, $\mathbf{ki}$, $\mathbf{ij}$, $\mathbf{jk}$, ...), it has
respective coefficients ($0,dx,0,0,0,0,dxdy,0,0,...$). And, as a member of
the structure spanned by ($1,dx$, $dy$, $dz$, $dydz$, $dzdx$, $dxdy$, ...),
it has coefficients ($0,$ $\mathbf{i}$, $0$, $0$, $0$, $0$, $\mathbf{ij}$, $%
0,0,\mathbf{...}).$ This is to be compared with, say, $dx\mathbf{j}$, whose
components in one algebra are (($0,0,dx,0,0,0,...)$ and ($0,$ $\mathbf{j}$, $%
0$, $0$, $0$, $0...)$ in the other (third place versus second place).

What is the motivation for considering the structure of mirror elements,
either in itself or as a substructure of the aforementioned tensor product
of structures? If the significance of $dx\mathbf{i+}dy\mathbf{j+}dz\mathbf{k}
$ were not enough to answer this question, consider rotations and associated
gauge transformations. The unit imaginary $i$ in the exponent in standard
treatments usurps the role of real geometric units of square minus one in
the tangent Clifford algebra of $E_{3}$. Thus ``$i$'' unduly replaces $%
\mathbf{a}_{1}\wedge \mathbf{a}_{2}$ ($=\mathbf{a}_{1}\mathbf{a}_{2}$ by
virtue of orthogonality) in the expressions for rotations around the $z$
axis when one uses Clifford algebra \cite{Casanova}, \cite{Lounesto}, \cite%
{Snygg}. Consequently, the replacement of $i$ with $\mathbf{a}_{1}\mathbf{a}%
_{2}$ should also take place in the phase factor for third component of
angular momentum and in the idempotents for solutions with rotational
symmetry. Under such a replacement, elements like $dxdy\mathbf{a}_{1}\mathbf{%
a}_{2}$ (no sum) emerge in the rotational idempotents.

Let $\Omega ^{lm}$ be differential $2$-forms, say the components $%
R_{pq}^{lm}dx^{p}\wedge dx^{q}$ of the curvature of Euclidean connections as
bivectors. The $\Omega ^{lm}$ are not mirror-symmetric, but the $%
R_{pq}^{pq}dx^{p}\wedge dx^{q}\mathbf{e}_{p}\mathbf{\wedge e}_{q}$ are, both
with and without summation over $p$ and $q$ (We use $\mathbf{e}$ instead of $%
\mathbf{a}$ to signify that the bases may be different at different points).

The set of mirror symmetric elements is closed under mirror symmetric
products, i.e. with the same symbol in both algebras. $\mathbf{(\wedge
,\wedge )}$ and $\mathbf{(\cdot ,\cdot )}$ are mirror products, but not $%
\mathbf{(\vee ,\wedge ),(\vee ,\cdot ),(\wedge ,\cdot ),(\wedge ,\vee ),}$%
etc. The little computation 
\begin{equation}
dy\mathbf{j(\vee ,\vee )}dx\mathbf{\mathbf{i}=}dydx\mathbf{ji=(-}dxdy)(-%
\mathbf{ij})=dxdy\mathbf{ij}  \label{57}
\end{equation}%
exhibits both the meaning of the symbol $\mathbf{(\vee ,\vee )}$ and the
commutativity of this product. This commutativity has nothing to do with the
fact that we multiplied monomials. Thus, for example, we have%
\begin{equation*}
(dx\mathbf{i+}dxdy\mathbf{\mathbf{ij}})\mathbf{(\vee ,\vee )}dxdz\mathbf{%
ik\equiv }dxdxdz\mathbf{\mathbf{iik+}}dxdydxdz\mathbf{\mathbf{ij}ik=}
\end{equation*}%
\begin{equation}
\mathbf{=}dz\mathbf{\mathbf{k+(-}}dydz)(-\mathbf{\mathbf{j}k})=dz\mathbf{%
\mathbf{k+}}dydz\mathbf{\mathbf{j}k=}dxdz\mathbf{ik(\vee ,\vee )}(dx\mathbf{%
i+}dxdy\mathbf{\mathbf{ij}}).
\end{equation}%
Coefficients that are not scalar functions do not change the mirror symmetry
element property. Thus%
\begin{equation*}
(\mu dx\mathbf{i+}\nu dxdy\mathbf{\mathbf{ij}})\mathbf{(\vee ,\vee )}\lambda
dxdz\mathbf{ik}=\mu \lambda dxdxdz\mathbf{\mathbf{iik+}}\nu \lambda dxdydxdz%
\mathbf{\mathbf{ij}ik=}
\end{equation*}%
\begin{equation}
\mathbf{=}\mu \lambda dz\mathbf{\mathbf{k}}+\nu \lambda \mathbf{\mathbf{(-}}%
dydz)(-\mathbf{\mathbf{j}k})=\mu \lambda dz\mathbf{\mathbf{k+}}\text{ }\nu
\lambda dydz\mathbf{\mathbf{j}k}.
\end{equation}%
We have just illustrated that, under mirror-symmetric products,
mirror-symmetric elements constitute a commutative algebra. It is a
restriction of the aforementioned tensor product of Clifford algebras. The
reason is that any minus signs emerging in a product of monomials in one
algebra is cancelled by a minus sign from the product of mirror images. The
non-symmetric product of mirror symmetric elements is not a mirror symmetric
element.

We shall refer with the names mirror symmetric algebra, commutative algebra
and commutative subalgebra to the restriction of that tensor product of
Clifford algebras to the set of mirror symmetric elements endowed with
mirror symmetric multiplication.

Given the element $dx\mathbf{+}dxdy$ in the K\"{a}hler algebra, it uniquely
determines the mirror element $dx\mathbf{i+}dxdy\mathbf{\mathbf{ij}}$. For
comparison of meanings, we observe that $\mathbf{a}_{1}\mathbf{a}_{2}$
represents any geometric object of unit area in the $xy$ plane of any given
Euclidean space. $dxdy$ represents an integrand whose integration gives a
number, namely the area of a figure contained in the $xy$ plane. And $dxdy%
\mathbf{a}_{1}\mathbf{a}_{2}$ represents an integrand whose integration
gives a bivector $\lambda \mathbf{a}_{1}\mathbf{a}_{2}$, namely any
geometric object of area $\lambda $ in the $xy$ plane of any given Euclidean
space with a reference frame attached to it.

\subsection{On the foundations of Clifford algebra.}

We have developed Clifford-related structure that has not been exploited by
the community of Clifford practitioners. One fundamental idea in this paper
is the replacement of tangent Clifford algebra and of K\"{a}hler algebra,
whether scalar or tensor-valued, by a tangent Clifford algebra valued
differential forms. The latter is the tensor product of the tangent Clifford
algebra by the Clifford algebra of scalar-valued differential forms. A
second fundamental idea is the realization that a commutative algebra lies
hidden in that structure when one just considers what we have called mirror
symmetric elements, and one restricts their multiplications to
mirror-symmetric multiplications. Because of the significance of the
implications of these fundamental ideas ---product of Clifford structures
and multiplication of privileged pairs--- they should be considered as
pertaining to the foundations of Clifford algebra.

Those two ideas can already be found in Helmstetter \cite{Helms}, though he
does not make a claim of having fundamental character. He, however, deals
with infinite dimensional vector spaces and his topics\ (Clifford monoids,
Lipschitz monoids and groups and deformations of Clifford algebras) have
nothing in common with ours.

In his development of the concept of Clifford monoid, Helmstetter gives the
multiplication%
\begin{equation}
(x\otimes y)\cdot (x^{\prime }\otimes y^{\prime })=\sigma (y,x^{\prime
})\sigma (y,y^{\prime })\cdot xx^{\prime }\otimes y^{\prime }y,  \label{60}
\end{equation}%
where the $\sigma $'s are just signs ($\pm )$ that depend on the arguments
(contents of the parentheses) and where ($x,y$) and ($x^{\prime },y^{\prime
} $) are members of the tensor product of two copies of the same Clifford
algebra. $x$ and $y$ belong to respective copies. But this is a by-product;
the pairs ($x,y$) and ($x^{\prime },y^{\prime }$) to which Eq. (60) refers
start as privileged pairs of elements of $A$. The tensor product on the
right hand side of (60) is not an hoc multiplication, but implied by an
abstract development whose underlying motivation will escape those who, like
the present author, lack the required algebraic sophistication. Without need
to understand the intricacies of this equation and of the argument that led
to it, one may nevertheless see similarities and differences. We mean the
likes of the structural similarity of (60) to our tensor product of
algebras. until we hit differences like Helmstetter's $xx^{\prime }\otimes
y^{\prime }y$ versus our own $xx^{\prime }\otimes yy^{\prime }.$ Beyond
that, we would be dealing with differences.

There is a variety of reasons, both mathematical and physical, that have led
us in past publications to the structures introduced in this paper. We do
not wish to enter them here, as there is reason enough for them in this
paper. However, approaching something as meaningful as a commutative algebra
from non-commutative ones through tensor product and then restricting both
the elements and the products may lack a depth which would perhaps be
achieved by other means. Just for the sake of the example as to the
existence of alternative approaches, but without making any claims as to
fundamental significance, one could start with elements like ($dx\mathbf{i}$%
, $dy\mathbf{j}$ and $dz\mathbf{k}$). One would them simply multiply them
with a Clifford product. Because of how the structure of what we have called
commutative algebra is born (the $dx\mathbf{i}$, $dy\mathbf{j}$ and $dz%
\mathbf{k}$ being indivisible elements rather than products involving $dx$
and $\mathbf{i}$), there would not be any need to speak of mirror symmetric
products of pre-existing elements. If there is fundamental substance to
these remarks on our construction, one should assume that there will also be
fundamental substance to Helmstetter's work, given the similarities
mentioned.

Let us mention another fundamental idea. Recall that a Euclidean space is a
special case of affine space, which in turn is a particular case of
projective space, in turn of ..., in turn of topological\textbf{\ }space, in
turn of ... and in turn of set (Fill the blanks as per your understanding).
Ideally, algebraic concepts at one of those levels (say Clifford algebra at
the Euclidean and pseudo-Euclidean levels) should mimic algebraic concepts
of more general levels. This point has been made very explicitly by
Barnabei, Brini and the late Professor Rota \cite{Rota} in connection with
the work of Peano \cite{Peano} on the algebra which is to projective space
what Clifford algebra is to Euclidean space. Grossly speaking, they make the
point that basic products in projective algebra should mimic operations in
set theory like those of reunion and intersection. Such conceptual
prolongation is already present, though in less overreaching way, in \'{E}.
Cartan's view of the theory of Euclidean connections as theory of affine
connections where affine tangent bundles are restricted to Euclidean tangent
bundles \cite{Cartan23}.

There is in Barnabei et al's work the attempt to look at mathematics from a
perspective of unification of conceptual spaces. We find such attempts very
laudable. At the level of trying to find such unifications, or at least
mentioning their potential existence, the issue of whether we are dealing
with finite or infinite spaces will certainly have consequences, but these
may come in the wash of developing the hypothetical unified picture.

\subsection{Action of angular momentum on the mirror algebra of 3-D
Euclidean space}

Of special significance is the action of angular momentum on rotational
idempotents in the mirror algebra.

Define%
\begin{equation}
\mathbf{I}_{ij}^{\pm }\equiv \pm \frac{1}{2}(1\pm w_{k}\mathbf{a}_{ij}).
\label{61}
\end{equation}%
Instead of (33) and (35) we now have%
\begin{equation}
J_{l}\mathbf{I}_{ij}^{\pm }=\pm \frac{1}{2}J_{l}w_{k}\mathbf{a}_{ij},
\end{equation}%
where $\mathbf{a}_{ij}$ is compact notation for $\mathbf{a}_{i}\wedge 
\mathbf{a}_{j}$ (equal to $\mathbf{a}_{i}\mathbf{a}_{j}$ by virtue of
orthonormality). Hence%
\begin{equation}
J_{i}\mathbf{I}_{ij}^{\pm }=\pm \frac{1}{2}w_{j}\mathbf{a}_{ij},\text{ \ \ \
\ }J_{j}\mathbf{I}_{ij}^{\pm }=\mp \frac{1}{2}w_{i}\mathbf{a}_{ij},\text{ \
\ \ \ }J_{k}\mathbf{I}_{ij}^{\pm }=0,
\end{equation}%
and, therefore,%
\begin{equation}
(K+1)\mathbf{I}_{ij}^{\pm }=\pm \frac{1}{2}(K+1)w_{k}\mathbf{a}_{ij},
\end{equation}%
since the action of $(K+1)$ on $1/2$ is zero.

We proceed to develop $(K+1)\mathbf{I}_{12}^{\pm }$ and generalize by means
of cyclic permutations. Thus%
\begin{equation}
(K+1)\mathbf{I}_{12}^{\pm }=\pm \frac{1}{2}\sum_{l}(J_{l}w_{3}\mathbf{a}%
_{12})w_{l}=\pm \frac{1}{2}(w_{2}w_{1}-w_{1}w_{2})\mathbf{a}_{12}=\pm w_{3}%
\mathbf{a}_{12},  \label{65}
\end{equation}%
and, therefore,%
\begin{equation}
(K+1)\mathbf{I}_{ij}^{\pm }=\pm w_{k}\mathbf{a}_{ij}.
\end{equation}%
The action of $K$ itself on $\mathbf{I}_{ij}^{\pm }$ is of interest:%
\begin{equation}
K\mathbf{I}_{ij}^{\pm }=\pm w_{k}\mathbf{a}_{ij}-\frac{1}{2}(1\pm w_{k}%
\mathbf{a}_{ij})=-\frac{1}{2}\pm \frac{1}{2}w_{k}\mathbf{a}_{ij}=-\mathbf{I}%
_{ij}^{\mp }.
\end{equation}%
Hence,%
\begin{equation}
-K\mathbf{I}_{ij}^{\pm }=\mathbf{I}_{ij}^{\mp },\text{ \ \ \ \ \ \ \ \ \ }%
K^{2}\mathbf{I}_{ij}^{\pm }=\mathbf{I}_{ij}^{\pm }.
\end{equation}%
We see that the rotational idempotents, $\mathbf{I}_{ij}^{\pm }$, are not
proper functions of either $K$ or $K+1$. The $w_{k}\mathbf{a}_{ij}$ are:%
\begin{equation}
(K+1)w_{k}\mathbf{a}_{ij}=2w_{k}\mathbf{a}_{ij},\text{ \ \ \ \ \ \ \ \ }%
Kw_{k}\mathbf{a}_{ij}=w_{k}\mathbf{a}_{ij}.
\end{equation}

For later purposes, we proceed to consider proper values of elements of
different grades in a Cartesian basis of differential $1-$forms and $2-$%
forms. Let $d\mathbf{x}^{i}$ denote $dx^{i}\mathbf{a}_{i}$ (no sum). We have%
\begin{equation}
J_{i}d\mathbf{x}^{1}=(0,\text{ }dx^{3},-dx^{2})\mathbf{a}_{1}  \label{70}
\end{equation}%
and 
\begin{equation}
(K+1)d\mathbf{x}^{1}=\sum (J_{i}d\mathbf{x}%
^{1})w^{i}=(dx^{3}w_{2}-dx^{2}w_{3})\mathbf{a}_{1}=2dx^{1}\mathbf{a}_{1}=2d%
\mathbf{x}^{1},
\end{equation}%
and similarly for other values of the index. Hence%
\begin{equation}
Kd\mathbf{x}^{l}=d\mathbf{x}^{l}.
\end{equation}%
On the other hand%
\begin{equation}
(K+1)dx^{ij}\mathbf{a}_{ij}=2dx^{ij}\mathbf{a}_{ij},
\end{equation}%
and, therefore,%
\begin{equation}
Kdx^{ij}\mathbf{a}_{ij}=dx^{ij}\mathbf{a}_{ij}.  \label{74}
\end{equation}%
It follows that the proper values of $(K+1)$ are always $+1$ for all members
in the module of mirror $1-$forms and mirror $2-$forms. They are zero for $%
0- $forms and mirror $3-$forms (Like scalars, $w$ commutes with the whole
algebra).\ The action of $K+1$ on it is, therefore, zero.

\section{Idempotents in Commutative Algebra}

\subsection{Hierarchy of products of idempotents in commutative algebra}

In the previous section, we introduced the three pairs of rotational
idempotents $\mathbf{I}_{ij}^{\pm }$ (i.e. $\pm \frac{1}{2}(1\pm w_{k}%
\mathbf{a}_{ij})$). Although commutativity is no longer a problem, once we
have one of those pairs there is no room for the other two in a product.

In a way similar to our definition of rotational idempotents in the
commutative algebra, we now define time and space translation symmetries.
Instead of the first equations 3, we now have, for the first of those two,%
\begin{equation}
\mathbf{\varepsilon }^{\pm }=\frac{1}{2}(1\mp dt\mathbf{e}_{0})\equiv \frac{1%
}{2}(1\mp d\mathbf{t}),
\end{equation}%
with $d\mathbf{t\equiv }dt\mathbf{e}_{0}.$

Space-translational symmetry is different, in that the phase factor is not
of the standard type, since the square of $\mathbf{a}_{i}$ is not minus one.
For this reason, time translational symmetry has pre-eminence over it.

Commutativity enters the concept of primitiveness. Since it is present
everywhere in this algebra, it is no longer a significant concept. The
replacement $\frac{1}{2}(1\pm dz\mathbf{a}_{3})$ for the idempotents $\frac{1%
}{2}(1\pm idz)$ will now have a role to play. As we did with the rotational
idempotents, we shall consider all three directions,%
\begin{equation}
\mathbf{P}_{i}^{\pm }\equiv \frac{1}{2}(1\pm d\mathbf{x}^{i}),
\end{equation}%
with $d\mathbf{x}^{i}\equiv dx^{i}\mathbf{a}_{i}.$

All these idempotents will have to be accompanied by corresponding geometric
phase factors, 
\begin{equation}
e^{m\phi \mathbf{a}_{i}\mathbf{a}_{j}}\mathbf{I}_{ij}^{\pm },\text{ \ \ \ \
\ }e^{-E\tau \mathbf{a}_{0}}\mathbf{\varepsilon }^{\pm }\text{\ \ \ \ \ \ }%
e^{\lambda _{i}x^{i}\mathbf{a}_{i}}\mathbf{P}_{i}^{\pm },\text{ \ \ \ \ \ no
sums!}  \label{77}
\end{equation}

We might think of proceeding with $d\rho $ as we have done with $dz.$ We do
not see a reason at this point to do so in dealing with ternary idempotents.
From a physical perspective, ``$d\rho $ translational symmetry'' is unheard
of. From a mathematical perspective, $dz$ and the corresponding idempotent
are constant differentials; $d\rho $ and its associated idempotents are not.
Also, the metric of Euclidean space depends on $\rho $ when written in
cylindrical coordinates (We have to use these coordinates if rotational
symmetry is assumed in the first place).

Needless to say that $\mathbf{I}_{ij}^{\pm }\mathbf{P}_{k}^{\pm }$ may be
expected to play a more significant role than $\mathbf{I}_{ij}^{\pm }\mathbf{%
P}_{i}^{\pm }$ and $\mathbf{I}_{ij}^{\pm }\mathbf{P}_{j}^{\pm }.$ Similarly,
in dealing with applications, products involving $\mathbf{I}_{ij}^{\pm }$'s
and/or $\mathbf{\varepsilon }^{\pm }$'s will be more significant than
products involving $\mathbf{P}_{i}^{\pm }$'s because of the inadequacy of
the corresponding exponential factors as phase factors.

Hence, we have been able to justify consideration of a set of idempotents
for symmetry which goes far beyond the set of K\"{a}hler's four binary
idempotents for dealing with the electromagnetic interaction.

We have not yet considered boosts. There are physical reasons for not doing
so, but this is not the place to deal with such reasons. Let us deal with
the mathematical arguments. Recall that, in dimension four, we had at our
disposal $\rho $ and $z$ ---but not on $\phi $ and $t$--- once joined time
translational and rotational symmetry had been assumed. There is no room to
accommodate any boost under such restrictions. One might retort that there
is no room either for translations along the $x$ and $y$ axes. But there is
an important distinction between boosts and translations, the implications
of this difference to be considered in a future paper. Let us give a hint.
Although there is no room for solutions of field equations with
translational symmetry in directions contained in the plane of rotations,
there still is room for them in connection with particles, specifically when
describing curves. Indeed curves satisfy the ``natural lifting conditions'' $%
dx^{i}-u^{i}dt=0$. So, $dx^{i}$ is just a multiple of $dt$ (or of $d\tau $).
Nothing similar takes place with the differentials $dx^{i}dt$ that pertain
to boosts like the $dx^{i}dx^{j}$'s pertain to rotations.

We are thus led to view the ternary idempotents%
\begin{equation}
\mathbf{\varepsilon }^{\pm }\mathbf{I}_{ij}^{\pm }\mathbf{P}_{m}^{\pm }
\end{equation}%
as candidates for physical applications in the same way as K\"{a}hler
considered $\epsilon ^{\pm }I_{xy}^{\ast }$ for applications in relativistic
quantum mechanics. We shall refer to (78) as the canonical ternary
idempotents of 4-D Lorentzian space.

\subsection{Canonical ternary idempotents of 4-D Lorentzian space}

In connection with (78), there are two issues to be considered. One of them
is that some of those are repeated. The other one is that, after removing
repetitions, there still are linear combinations of them that add to zero.
We shall deal only with the first of these issues, since dealing with the
second one would have too much of a physics flavor. That is purposefully
avoided to avoid conflicts of classification in the arXiv.

Consider the following equalities:%
\begin{equation}
4\mathbf{I}_{12}^{+}\mathbf{P}_{1}^{\pm }=(1+d\mathbf{x}d\mathbf{y})(1\pm d%
\mathbf{x})=(1+d\mathbf{x}d\mathbf{y})(1\pm d\mathbf{y})=4\mathbf{I}_{12}^{+}%
\mathbf{P}_{2}^{\pm }.
\end{equation}%
\begin{equation}
4\mathbf{I}_{12}^{-}\mathbf{P}_{1}^{\pm }=(1-d\mathbf{x}d\mathbf{y})(1\pm d%
\mathbf{x})=(1-d\mathbf{x}d\mathbf{y})(1\mp d\mathbf{y})=4\mathbf{I}_{12}^{-}%
\mathbf{P}_{2}^{\mp }.  \label{80}
\end{equation}%
Notice the change of sign in the exponent of (80), which is not the case in
the exponent of (79). At the level of solutions of equations and not just
idempotents, there would still be a difference in the phase factors.

Multiplication of $\mathbf{I}_{12}$ with $\mathbf{P}_{3}^{\pm }$ will not be
in the same footing with multiplication with $\mathbf{P}_{1}^{\pm }$ and $%
\mathbf{P}_{2}^{\pm }$, since $z$ (thus $\mathbf{P}_{3}^{\pm }$) relates to $%
\mathbf{I}_{12}$ in a different way than to $\mathbf{P}_{1}^{\pm }$ and $%
\mathbf{P}_{2}^{\pm }$. In the case of $\mathbf{P}_{3}$, we have a different
symmetry if the plus and minus directions are indistinguishable. In that
case, $\mathbf{I}_{12}^{+}\mathbf{P}_{3}^{+}$ is equivalent to $\mathbf{I}%
_{12}^{-}\mathbf{P}_{3}^{-}$, and $\mathbf{I}_{12}^{+}\mathbf{P}_{3}^{-}$ is
equivalent to $\mathbf{I}_{12}^{-}\mathbf{P}_{3}^{+}$ (There are subtleties
involved in this which should be taken into account in applications where
more than idempotent would be involved. Again, we stay out of this since we
would be entering physics.

In view of these symmetries, we can organize the six different $\mathbf{I}%
_{12}^{{}}\mathbf{P}$ binary idempotents in three columns and two rows
respectively for the subscripts and superscripts of $\mathbf{P}$. We
left-multiply them by $\mathbf{\varepsilon }^{\pm }.$ We thus have the
following table of twelve different canonical ternary idempotents%
\begin{equation}
\left[ 
\begin{array}{ccc}
\text{ \ \ }\mathbf{\varepsilon }^{\pm }\mathbf{I}_{12}^{+}\mathbf{P}_{1}^{+}%
\text{ \ \ } & \text{ \ \ }\mathbf{\varepsilon }^{\pm }\mathbf{I}_{12}^{-}%
\mathbf{P}_{2}^{+}\text{ \ \ } & \text{ \ \ }\mathbf{\varepsilon }^{\pm }%
\mathbf{I}_{12}^{+}\mathbf{P}_{3}^{+}\text{ \ \ } \\ 
\text{ \ \ }\mathbf{\varepsilon }^{\pm }\mathbf{I}_{12}^{+}\mathbf{P}_{1}^{-}%
\text{ \ \ } & \mathbf{\varepsilon }^{\pm }\mathbf{I}_{12}^{-}\mathbf{P}%
_{2}^{-} & \mathbf{\varepsilon }^{\pm }\mathbf{I}_{12}^{+}\mathbf{P}_{3}^{-}%
\end{array}%
\right] .
\end{equation}%
We could also have organized the same group of idempotents with $\mathbf{I}%
_{12}^{-}$ in columns 1 and 3, and $\mathbf{I}_{12}^{+}$ in column 2,
everything else remaining the same. Other options would involve physical
considerations.

What we have done for $\mathbf{I}_{12}$ can be repeated for $\mathbf{I}_{23}$
and $\mathbf{I}_{31}.$ In these cases, the role of $\mathbf{P}_{3}$ in (81)
would be played by $\mathbf{P}_{1}$ and $\mathbf{P}_{2}$ respectively.
Putting these options together yields a set of 36 different ternary
idempotents. They cannot be linearly independent since the dimension of the
commutative algebra is 16. We shall not enter into that since the results
would look very much like physics, and we intentionally want to let this
paper be of a purely mathematical nature.

\section{Concluding Remarks}

Among the many results obtained in this paper we mention those that
constitute its core. K\"{a}hler's idempotents for solutions with symmetry
contain the unit imaginary. We have reformulated them so that their now
geometric form allows for the removal of that unit. It happens in a
commutative algebra embedded in the tensor product of K\"{a}hler's real
algebra of scalar-valued differential forms with its dual Clifford algebra
of valuedness (i.e. tangent Clifford algebra). More importantly, K\"{a}%
hler's focus on solutions of equations rather than on the equations
themselves allows one to extend the manifestation of symmetry through a
greater set of more sophisticated idempotents.\bigskip 

ACKNOWLEDGEMENTS\bigskip

Conversations with Professors Oziewicz and Helmstetter are acknowledged.
Generous funding from PST\ Associates is deeply appreciated.

\end{document}